\documentclass[smallextended]{article}
\usepackage{amsmath,amssymb,bbm,calc,capt-of,ifthen}
\usepackage[english]{babel}
\usepackage{graphicx}
\usepackage{xcolor}
\usepackage{listings}
\usepackage{matlab-prettifier}

\begin{document}

\title{FEM and BEM simulations with the {\sc{Gypsilab}} framework}


\author{
Fran\c{c}ois Alouges\footnote{
              CMAP - Ecole Polytechnique, Universit\'e Paris-Saclay, Route de Saclay, 91128, Palaiseau Cedex, France. 			{\tt{francois.alouges@polytechnique.edu}}\,. } \,\,and 
            Matthieu Aussal\footnote{CMAP - Ecole Polytechnique, Universit\'e Paris-Saclay, Route de Saclay, 91128, Palaiseau Cedex, France. {\tt{matthieu.aussal@polytechnique.edu}\,.}}
}

\date{}

\maketitle

\begin{abstract}
{\sc{Gypsilab}} is a {\sc{Matlab}} framework which aims at simplifying the development of numerical methods
that apply to the resolution of problems in multiphysics, in particular, those involving FEM or BEM simulations. 
The peculiarities of the framework, with a focus on its ease of use, are shown together with the methodology 
that have been followed for its development. Example codes that are short though representative enough
are given both for FEM and BEM applications. A performance comparison with FreeFem++ is provided, and a
particular emphasis is made on problems in acoustics and electromagnetics solved using the BEM and for which
compressed $\mathcal{H}$-matrices are used.
\end{abstract}

\section{Introduction}
The finite element method (FEM) is nowadays extremely developed and has been widely used in the
numerical simulation of problems in continuum mechanics, both for linear and non-linear problems. 
Many software packages exist among which we may quote 
free ones (e.g. {\sc{FreeFem++}} \cite{FreeFem++}, FENICS \cite{Fenics}, FireDrake \cite{Firedrake}, 
Xlife++ \cite{Xlife++}, Feel++ \cite{Feel++}, GetDP \cite{getdp}, etc.)
or commercial ones (e.g. COMSOL \cite{comsol}). The preceding list is far from
being exhaustive as the method has known many developments and improvements and is still 
under active study and use. Numerically speaking, and without entering into too many details, let us say that the particularity of the method is that it is based on a weak formulation that leads to sparse matrices for which the footprint in memory is typically proportional to the number of unknowns. Direct or iterative methods can then be used to solve the underlying linear systems.

On the other hand, the boundary element method (BEM) is used for problems which can be expressed using a Green kernel.
A priori restricted to linear problems, the method inherently possesses the faculty of handling free space 
solutions and is therefore currently used for solving Laplace equations, wave scattering problems (in acoustics, electromagnetics or elastodynamics)
or Stokes equations for instance. Although it leads to dense matrices, which storage is proportional to the square
of the number of unknowns, the formulation is usually restricted to the boundary of the domain under consideration
(e.g. the scatterer), which lowers the dimension of the object that needs to be discretized. Nevertheless, due to computer limitations, those methods may require a special
compression technique such as the Fast Multipole Method (FMM) \cite{FMM,FMMpage}, the $\mathcal{H}$-matrix storage
\cite{hmatrix} or the more recent Sparse Cardinal Sine Decomposition \cite{SCSD1,SCSD2,SCSD3}, in order to be applied to relevant problems. In terms of software packages available for the simulation with this kind of methods, and probably due to the technicality sketched above, the list is much shorter than before. In the field of academic, freely available software packages, one can quote BEM++ \cite{bem++}, or Xlife++ \cite{Xlife++}. 
Commercial software packages using the method are for instance the ones distributed by ESI Group (VAone for acoustics \cite{vaone} and E-Field \cite{efield} for electromagnetism), or Aseris \cite{aseris}. Again, the preceding list is certainly not exhaustive.

Eventually, one can couple both methods, in particular to simulate multiphysics problems that involve different materials and for which none of the two methods apply directly. This increases again the complexity of the methodology as the user needs to solve coupled equations that are piecewise composed of matrices sparsely stored
combined with other terms that contain dense operators or compressed ones.  How to express such a problem? Which format should be used for the matrix storage? Which solver applies to such cases, an iterative or a direct one? Eventually, the writing of the software still requires abilities that might be out of the user's field of expertise.

The preceding considerations have led us to develop the framework {\sc{Gypsilab}}  which, in particular, aims at simplifying and generalizing the
development of FEM-BEM coupling simulation algorithms and methods. Written as a full {\sc{Matlab}} suite of routines, the
toolbox can be used to describe and solve various problems that involve FEM and/or BEM techniques. 
We have also tried to make the
finite element programming as simple as possible, using a syntax comparable to the one used in {\sc{FreeFem++}} 
or FENICS/FireDrake and very close to the mathematical formulation itself which is used to discretize the problems. 
Let us emphasize that, to our best knowledge, the only freely available software package, which uses the BEM and for which the way of programming is comparable to {\sc{FreeFem++}}, is BEM++ that has been interfaced into FENICS. The software described in this paper thus provide the user with a {\sc{Matlab}} alternative.  
Fairly easy to install since it is written in full {\sc{Matlab}} without any piece of C code,
the framework contains at the moment several toolboxes:
\begin{itemize}
\item {\sc{openMSH}} for the handling of meshes;
\item {\sc{openDOM}} for the manipulation of quadrature formulas and weak formulations that involve integrals;
\item {\sc{openFEM}} for the finite element and boundary element methods;
\item {\sc{openHMX}} that contains the full $\mathcal{H}$-matrix algebra \cite{hmatrix}, including the LU factorization and inversion.
\end{itemize}

In what follows, we will assume that all the preceding librairies are in the user's path. This could be done using the following commands, assuming that the librairies have been placed in the parent directory of the current one.
\vspace*{0.3cm}

\hrule
\begin{lstlisting}[style=Matlab-editor,basicstyle=\small]
% Library paths
addpath('../openMsh')
addpath('../openDom')
addpath('../openFem')
addpath('../openHmx')
\end{lstlisting}
\hrule
\vspace*{0.3cm}

Eventually, although the main goal is not the performance, {\sc{Gypsilab}} may handle, on classical computers, problems whose size reaches a few millions of unknowns for the FEM part and a few hundreds of thousands of unknowns for the BEM when one uses compression. For FEM applications, this performance is very much comparable to already proposed {\sc{Matlab}} strategies \cite{Cuvelier,Anjam,Rahman}, though {\sc{Gypsilab}} presents a much higher generality and flexibility. 

The present paper explains in some details the capabilities of {\sc{Gypsilab}} together with its use. Explanations concerning the implementation are also provided that give an insight into the genericity of the approach. In order to simplify the exposition, we focus here on FEM or BEM applications, leaving coupled problems to a forthcoming paper. 

Due to its simplicity of use, we strongly believe that the software library described here could become a reference in the field. Indeed, applied mathematicians, interested in developing new FEM-BEM coupled algorithms need such tools in order to address problems of a significant size for real applications. Moreover, the library is also ideal for the quick prototyping of academic or industrial applications.
  
\section{Simple examples}
We show in this Section a series of small example problems and corresponding {\sc{Matlab}} listings.
\subsection{A Laplace problem with Neumann boundary conditions}
Let us start with the writing of a finite element program that solves the partial differential equation
(PDE)
\begin{equation}
\left\{
\begin{array}{l}
-\Delta u +u =f \mbox{ on } \Omega\,,\\
\displaystyle \frac{\partial u}{\partial n} = 0 \mbox{ on }\partial \Omega\,,
\end{array}
\right.
\label{Laplace1}
\end{equation}
where the right-hand side function $f$ belongs to $L^2(\Omega)$. Here $\Omega$ stands for a domain in $\mathbb{R}^2$ or $\mathbb{R}^3$.

The variational formulation of this problem is very classical and reads:\\

Find $u\in H^1(\Omega)$ such that $\forall v\in H^1(\Omega)$
$$
\int_\Omega \nabla v(x) \cdot \nabla u(x)\,dx + \int_\Omega v(x) u(x)\,dx = \int_\Omega f(x) v(x)\,dx\,.
$$

The finite element discretization is also straightforward and requires solving the same variational formulation where $H^1(\Omega)$ is replaced by one of its finite dimensional subspaces (for instance the set of continuous and piecewise affine on a triangular mesh of $\Omega$ in the case of linear $P^1$ elements).

We give hereafter the {\sc{Gypsilab}} source code used to solve such a problem in the case where the domain under consideration is the unit disk in $\mathbb{R}^2$, and the function $f$ is given by $f(x,y)=x^2$.
\vspace*{0.3cm}

\hrule
\begin{lstlisting}[style=Matlab-editor,basicstyle=\small]
% Mesh of the disk
N = 1000;
mesh = mshDisk(N,1);
% Integration domain
Omega = dom(mesh,3);
% Finite elements
Vh = fem(mesh, 'P1');
% Matrix and RHS
f = @(X) X(:,1).^2;
K = integral(Omega,grad(Vh),grad(Vh)) + integral(Omega,Vh,Vh);
F = integral(Omega, Vh, f);
% Solving
uh = K \ F;
figure
graph(Vh,uh);
\end{lstlisting}
\hrule
\vspace*{0.3cm}

\begin{figure}[h]
\centerline{\includegraphics[width = 8cm]{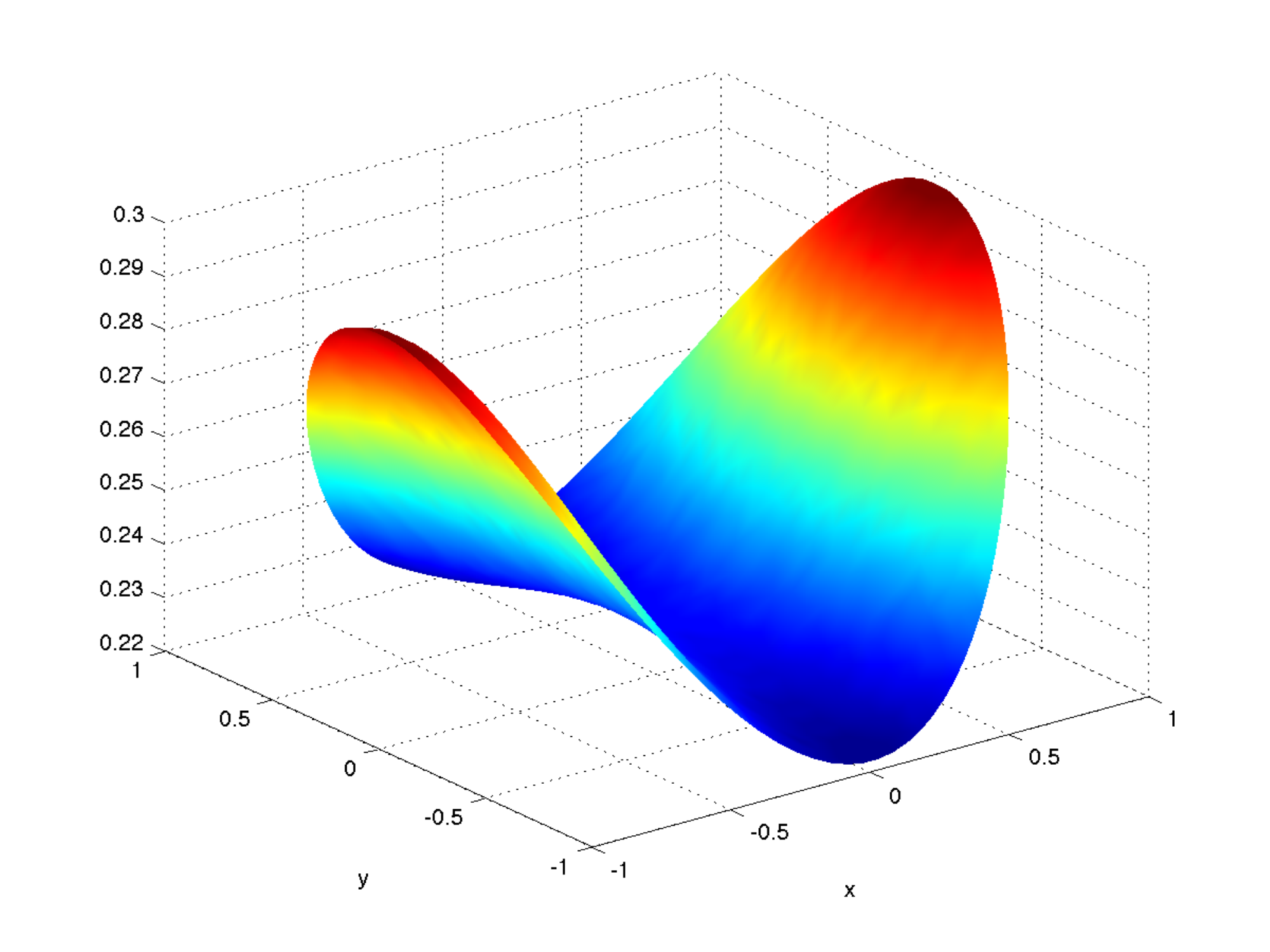}}
\caption{Numerical solution of (\ref{Laplace1}) on a unit disk using {\sc{Gypsilab}}.}
\label{fig1}
\end{figure}

We believe that the listing is very clear and almost self-explanatory. The disk is first meshed with 1000 vertices  (lines 2-3), then one defines an integration domain (line 5), the finite element space (line 7),
the weak formulation of the problem (lines 9-11) and the resolution (lines 13). 
Let us immediately insist on the fact that the operators constructed by the {\tt{integral}} keyword are 
really matrix and vector {\sc{Matlab}} objects so that 
one can use classical {\sc{Matlab}} functionalities for the resolution (here the ``backslash'' operator {\tt{\textbackslash}}). 
Plotting the solution (lines 14-15) leads to the figure reported in Fig. \ref{fig1}.

\subsection{Fourier boundary conditions}
We now consider the problem
\begin{equation}
\left\{
\begin{array}{l}
-\Delta u +u =0 \mbox{ on } \Omega\,,\\
\displaystyle \frac{\partial u}{\partial n} + u = g \mbox{ on }\partial \Omega\,,
\end{array}
\right.
\label{Laplace2}
\end{equation}

Again, the weak formulation of the problem is standard and reads as follows\\

Find $u\in H^1(\Omega)$ such that $\forall v\in H^1(\Omega)$
$$
\int_\Omega \nabla v(x) \cdot \nabla u(x)\,dx +\int_\Omega v(x)u(x)\,dx +\int_{\partial \Omega} v(s)u(s)\,ds= \int_{\partial \Omega} g(s) v(s)\,ds\,.
$$

The preceding code is modified in the following way (we have taken the example where $g(s)=1$).

\vspace*{0.3cm}

\hrule
\begin{lstlisting}[style=Matlab-editor,basicstyle=\small]
% Create mesh disk + boundary
N = 1000;
mesh = mshDisk(N,1);
meshb = mesh.bnd;
% Integration domains
Omega = dom(mesh,7);
Sigma = dom(meshb,3);
% Finite element space
Vh = fem(mesh,'P2');
% Matrix and RHS
K = integral(Omega, grad(Vh), grad(Vh)) ...
     + integral(Omega, Vh, Vh) ...
     + integral(Sigma, Vh, Vh);
g = @(x) ones(size(x,1),1);
F = integral(Sigma, Vh, g);
% Resolution
uh = K \ F;
\end{lstlisting}
\hrule
\vspace*{0.3cm}

\begin{figure}[h]
\centerline{\includegraphics[width = 9cm]{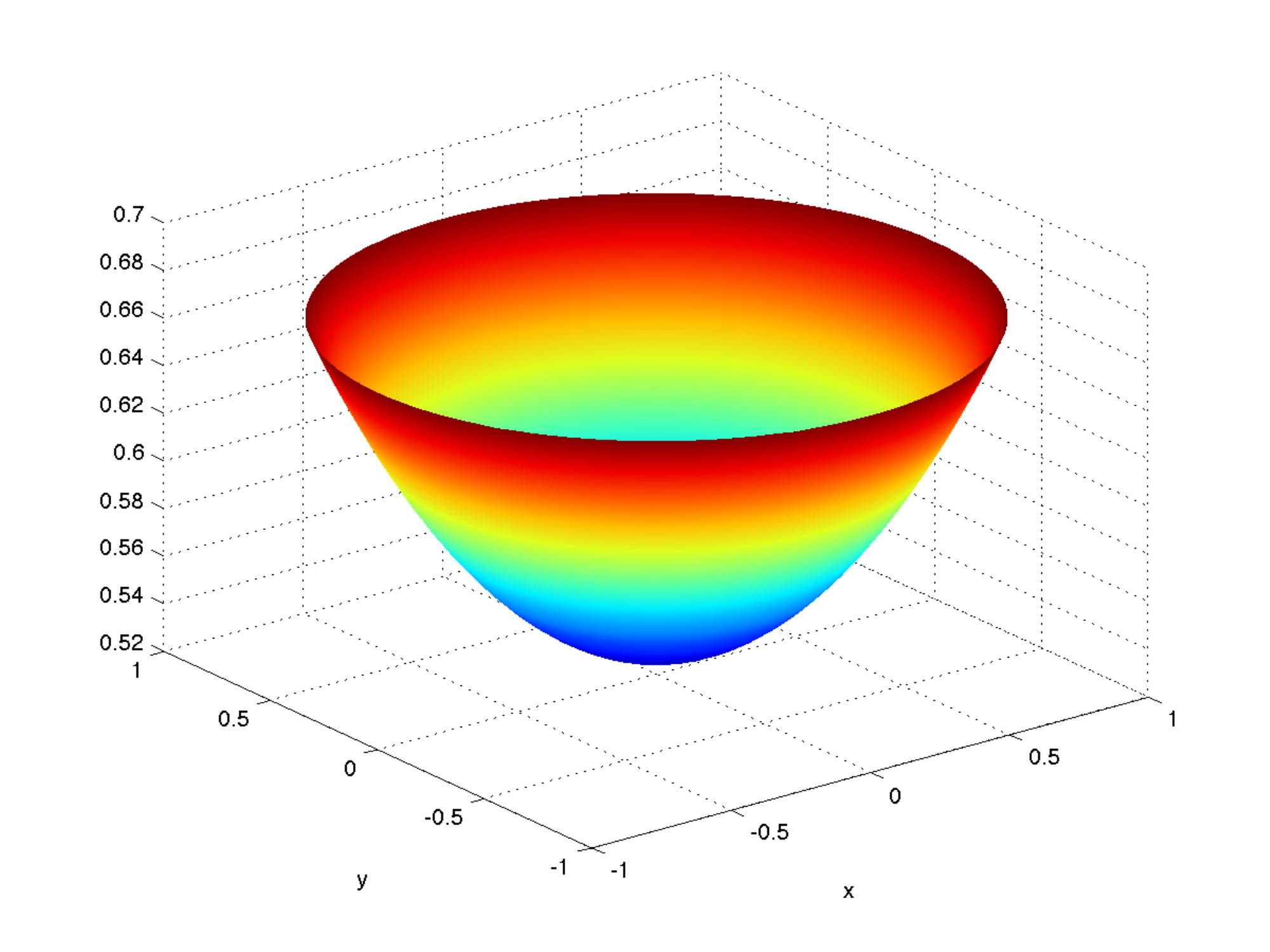}}
\caption{Numerical solution of (\ref{Laplace2}) on a unit disk using {\sc{Gypsilab}}.}
\label{fig2}
\end{figure}

Compared to the preceding example, a boundary mesh and an associated integration domain are also defined (lines 4 and 7). Let us note the piecewise quadratic (so-called $P^2$) element used (line 9) which leads to use more accurate integration formulas, respectively with 7 points for the triangles (line 6) and 3 points per segment for the boundary mesh (line 7). Again, the result obtained is plotted in Fig. \ref{fig2}.

\subsection{An eigenvalue problem}
We end up this section by showing an example of a 3D problem, namely, the computation of the first eigenvalues/eigenvectors
of the Laplace operator in the parallelepipedic volume $[0,1]\times[0,\frac12]\times[0,\frac12]$, with Dirichlet boundary conditions. Mathematically speaking the problem writes as finding couples $(\lambda, u)$ that are solutions to the eigenvalue problem
\begin{equation}
\left\{
\begin{array}{l}
-\Delta u = \lambda u \mbox{ on } \Omega\,,\\
\displaystyle u = 0 \mbox{ on }\partial \Omega\,,
\end{array}
\right.
\label{Laplaceeigen}
\end{equation}
where $\Omega=[0,1]\times[0,\frac12]\times[0,\frac12]$. Now the problem is posed in 3D, and we need to 
force (homogeneous) Dirichlet boundary conditions.
The corresponding {\sc{Matlab}} listing becomes
\vspace*{0.3cm}

\hrule
\begin{lstlisting}[style=Matlab-editor,basicstyle=\small]
% Create mesh of the cube + boundary
N = 1e4;
mesh = mshCube(N,[1 0.5 0.5]);
meshb = mesh.bnd;
% Integration domains
Omega = dom(mesh,4);
% Finite element space
Vh = fem(mesh,'P1');
Vh = dirichlet(Vh,meshb);
% Matrix and RHS
K = integral(Omega, grad(Vh), grad(Vh)) ;
M = integral(Omega, Vh, Vh);
% Resolution
Neig = 10;
[V,EV] = eigs(K,M,Neig,'SM');
\end{lstlisting}
\hrule
\vspace*{0.3cm}

Notice the enforcement of Dirichlet boundary conditions on the finite element space (line 9), the assembling of the rigidity and mass matrices (lines 11-12). The results obtained by the preceding code for the computation of the first 10 eigenelements of the Laplace operator with Dirichlet boundary condition is given in Table \ref{tableVP}.
\begin{table}[ht]
\begin{center}
\begin{tabular}{c|c|c|c}
\hline
Number & Exact & Numeric & Relative error\\
\hline
1 &      88.8264  & 90.1853   &  0.0153  \\
2 &    118.4353  & 120.5672 &  0.0180 \\
3 &    167.7833  & 171.4698 &  0.0220 \\
4 &    207.2617  & 213.3081 &  0.0292 \\
5 &    207.2617  & 213.4545 &  0.0299 \\
6 &    236.8705  & 243.2719 &  0.0270 \\
7 &    236.8705  & 244.4577 &  0.0320 \\
8 &    236.8705  & 245.0540 &  0.0345 \\
9 &    286.2185  & 296.6374 &  0.0364 \\
10 &  286.2185  & 297.9819 &  0.0411 \\
\hline
\end{tabular}
\end{center}
\caption{Exact and approximated eigenvalues of the Laplacian with Dirichlet boundary conditions on the parallelepipedic domain $[0,1]\times[0,\frac12]\times[0,\frac12]$. For each of the first ten eigenvalues, we give its exact value, the one computed with the program before and the relative error.} 
\label{tableVP}
\end{table}
\section{Finite element programming in {\sc{Matlab}}}
\label{sec3}
Programming the finite element method in {\sc{Matlab}} is very attractive and has already been
considered by many people (see for instance \cite{50lines,Funken,Sutton,Kwon,Cuvelier,Anjam,Rahman}), probably because the language is easy to use and already contains the most recent linear solvers.
It is important to notice that {\sc{Matlab}} is usually very powerful when one uses its vectorial capabilities. Therefore, the traditional assembling of the matrices that come from finite element discretizations, which uses a loop over the elements with indirect addressing, might lead to prohibitive execution times. This problem was identified a long time ago and
several ways have already been proposed to circumvent it. In particular, we refer the interested reader to  \cite{Cuvelier}, where different alternatives that lead to very efficient assembling are given and compared. Many languages are also compared (C++, {\sc{Matlab}}, {\sc{Python}}, {\sc{FreeFem++}}), and it is shown that the C++ implementation only brings a slight improvement in performance. Other {\sc{Matlab}} implementations are proposed in the literature (see e.g. \cite{50lines,Sutton,Kwon,Anjam,Rahman}), but they all suffer from the lack of generality. The problem solved is indeed very often the Laplacian with piecewise linear finite elements and one needs to adapt the approach for any different problem.

We have followed yet another strategy that has the great advantage to be very general and easily adaptable to a wide variety of possible operators to build and solve, and which also enables the user to assemble matrices that come from the coupling between different finite elements. Moreover, we will see that the method also leads to reasonably good assembling times. To this aim, we give the following example from which one can understand the generality of the approach and the way {\sc{Gypsilab}} is coded.

Let us consider the case of assembling the mass matrix. To be more precise, we call $\mathcal{T}$ a conformal triangulation\footnote{Triangulation usually means a 2D problem, while we would have to consider a tetrahedral mesh for 3D problems. This is not a restriction as we can see.} on which one has to compute the matrix $A$ whose entries are given by
 \begin{equation}
 A_{ij} = \int_{\mathcal{T}} \phi_i(x) \phi_j(x)\,dx\,.
 \label{mass}
 \end{equation}
 
 Here we have used the notation $(\phi_i)_{1\leq i\leq N}$ to denote the basis functions of the finite element (discrete) space
 of dimension $N$. To gain generality, the preceding integral is usually not computed exactly, but rather approximated using a quadrature rule. Thus, calling $(x_k,\omega_k)_{1\leq k\leq M}$ the set of all quadrature points over $\mathcal{T}$, we
 write
 \begin{equation}
 A_{ij} \sim \sum_{k=1}^M \omega_k \phi_i(x_k) \phi_j(x_k)\,.
 \label{approxA}
 \end{equation}
 
Introducing now the two matrices $W$ and $C$ (respectively of size $M\times M$ and $M\times N$) defined by
\begin{equation}
 W_{kk} = \omega_k \mbox{ for } 1\leq k\leq M\,,\mbox{ and } B_{kj} = \phi_j(x_k) \mbox{ for } 1\leq k\leq M,\,1\leq j\leq N\,,
\label{WB}
\end{equation}
we may rewrite (\ref{approxA}) as
\begin{equation}
A \sim B^T \,W\, B\,,
\label{A=BWB}
\end{equation} 
where $B^T$ denotes the transpose of $B$. Notice that the approximation is coming from the fact that a quadrature rule has been used instead of an exact formula. In particular,  we emphasize that if the quadrature formula is exact in (\ref{approxA}), then the approximation (\ref{A=BWB}) is in fact an equality.

From the preceding considerations the procedure that enables to assemble the sparse mass matrix can be summarized as:
\begin{itemize}
\item Knowing the triangulation (resp. tetrahedral mesh), and a quadrature formula on a reference triangle (resp. tetrahedron), build the set of quadrature points/weights $(x_k,\omega_k)_{1\leq k\leq M}$. (This is done in the package {\sc{openDOM}}.)
\item Knowing the finite element used (or, equivalently, the basis functions $(\phi_j)_{1\leq j\leq N}$) and the quadrature points $(x_k)_{1\leq k\leq M}$, build the matrices $W$ and $B$. (This is done in the package {\sc{openFEM}}.)
\item Eventually, compute $A = B^t \,W\,B\,$. 
\end{itemize}

Notice that the matrices $W$ and $B$ are usually sparse. Indeed, $W$ is actually diagonal, while $B$ has a non-zero entry $B_{jk}$ only for the quadrature points $x_k$ that belong to the support of $\phi_j$. In terms of practical implementation,
the matrix $W$ is assembled using a {\tt{spdiags}} command while the matrix $B$ is built using a vectorized 
technique. 

The preceding procedure is very general and does only rely on the chosen finite element or the quadrature formula. Moreover, the case of more complicated operators can also be treated with
only slight modifications. Indeed, if one considers the case of the Laplace operator, for which the so-called rigidity matrix is given by
 \begin{equation}
 K_{ij} = \int_{\mathcal{T}} \nabla \phi_i(x) \cdot \nabla \phi_j(x)\,dx\,,
 \label{rigidity}
 \end{equation}
one may write similarly
 \begin{equation}
 K_{ij} \sim \sum_{k=1}^M \omega_k \nabla \phi_i(x_k)\cdot  \nabla \phi_j(x_k)\,,
 \label{approxK}
 \end{equation}
from which one deduces
$$
K \sim C^t_x\,W\,C_x + C^t_y\,W\,C_y + C^t_z\,W\,C_z\,,  
$$
where the matrix $W$ is the same as before and the matrices $C_x,C_y$ and $C_z$ are given for $1\leq k\leq M,\,1\leq j\leq N$ by
\begin{equation}
 C_{x,kj} = \frac{\partial\phi_j}{\partial x}(x_k)\,,\,\,C_{y,kj} = \frac{\partial\phi_j}{\partial y}(x_k)\,,\,\,C_{z,kj} = \frac{\partial\phi_j}{\partial z}(x_k)\,.
\label{CxCyCz}
\end{equation}

\section{Quick overview of {\sc{Gypsilab}}}
This section is not intended to be a user's manual. We just give the main functionalities of {\sc{Gypsilab}}  and 
refer the interested reader to the website \cite{Gypsilab}. It is important to notice that {\sc{Gypsilab}} tries to compute as much as possible the quantities ``on the fly'', or in other words, to
keep in memory as little information as possible. The main underlying idea is that storing many vectors and matrices (even sparse) might become memory consuming, and recomputing on demand the corresponding quantities does not turn to be the most costly part in usual computations. Keeping this idea in mind helps in understanding the ``philosophy'' that we have followed for the development of the different toolboxes. Moreover, the whole library is object oriented and the toolboxes have been implemented as value classes.
\subsection{The mesh}
Finite or Boundary element methods are based on the use of a mesh. The routines that handle the mesh object are grouped into the toolbox {\sc{openMSH}}. In {\sc{Gypsilab}}, the mesh is a purely geometric object with which one can compute only geometric quantities (e.g. normals, volumes, edges, faces, etc.). A mesh can be of dimension 1 (a curve), 2 (a surface) or 3 (a volume) but is always embedded in the dimension 3 space and is a simplicial mesh (i.e. composed of segments, triangles or tetrahedra). It is defined by three tables~:
\begin{itemize}
\item A list of vertices, which is a table of size $N_v\times 3$ containing the three-dimensional coordinates of the vertices~;
\item A list of elements, which is a table of size $N_e\times (d+1)$, $d$ being the dimension of the mesh and $N_e$ the number of elements~;
\item A list of colors, which is a column vector of size $N_e\times 1$ defining a color for each element of the mesh, this last table being optional.
\end{itemize}  

A typical {\tt{msh}} object is given by the following structure.

\begin{verbatim}
>>mesh

mesh = 

  2050x3 msh array with properties:

    vtx: [1083x3 double]
    elt: [2050x3 double]
    col: [2050x1 double]
\end{verbatim}
    
The {\sc{openMSH}} toolbox does not yet contain a general mesh generator per se. Only simple objects (cube, square, disk, sphere, etc.) can be meshed. More general objects may be nevertheless loaded using classical formats ({\tt{.ply}}, {\tt{.msh}}, {\tt{.vtk}}, {\tt{.stl}}, etc.). Since the expected structure for a mesh is very simple, the user may also his/her own wrapper to create the previous tables.

The {\sc{openMSH}} toolbox also contains many operations on meshes such as the intersection or the union of different meshes, the extraction of the boundary, etc. Let us emphasize that upon loading, meshes are automatically cleaned by removing unnecessary or redundant vertices or elements.

\subsection{The domain}
The {\bf{domain}} is a geometric object on which one can furthermore integrate. Numerically speaking, this is the concatenation of a mesh and a quadrature formula. This formula is identified by a number, that one uses with the simplices of the corresponding mesh, in order to integrate functions. The default choice is a quadrature formula with only one integration point located at the center of mass of the simplices. This is usually very inaccurate, and it is almost always mandatory to enhance the integration by taking a higher degree quadrature formula. A domain is defined using the {\tt{dom}} keyword. For instance, the command
\begin{lstlisting}[style=Matlab-editor]
Omega = dom(myMesh,4);
\end{lstlisting}
defines an integration domain {\tt{Omega}} from the mesh {\tt{myMesh}}, using an integration formula with 4 integration points. If such an integration formula is not available, the program returns an error. Otherwise, the command creates an integration domain
with the structure shown by the following output.
\begin{verbatim}
>> omega = dom(myMesh,4)

omega = 

  dom with properties:

    gss: 4
    msh: [2050x3 msh]
\end{verbatim}

We believe that making the construction of the quadrature formula very explicit helps the user to pay attention to this very important point, and make the right choice for his/her application. In particular, for finite element computing, the right quadrature formula that one needs to use depends on the order of the chosen finite element. Integration functionalities are implemented in the {\sc{openDOM}} toolbox (see more below).

\subsection{The Finite Element toolbox ({\sc{openFEM}})}
Finite element spaces are defined through the use of the class constructor {\tt{fem}}. Namely, the command
\begin{lstlisting}[style=Matlab-editor]
Vh = fem(mesh, name);
\end{lstlisting}
creates a finite element space on the {\tt{mesh}} (an arbitrary 2D or 3D mesh defined in $\mathbb{R}^3$) of type defined by {\tt{name}}. At the moment of the writing of this paper, 3 different families of finite elements are available:
\begin{itemize}
\item The Lagrange finite elements. The orders 0, 1 and 2 are only available for the moment. They correspond to piecewise polynomials of degree 0, 1 and 2 respectively. 
\item The edge N\'ed\'elec finite element. It is a space of vectorial functions whose degrees of freedom are the circulation along all the edges of the underlying mesh. This finite element is defined in both 2D and 3D. In 2D is implemented a general form for which the underlying surface does not need to be flat.
\item The Raviart-Thomas, also called Rao-Wilton-Glisson (RWG) finite elements. Also vectorial, the degrees of freedom are the fluxes through the edges (in 2D) or the faces (in 3D) of the mesh. Again, the 2D implementation is available for general non-flat surfaces.
\end{itemize}
For the two last families, only the lowest orders are available. The value of the variable {\tt{name}} should be one of {\tt{'P0'}}, {\tt{'P1'}}, {\tt{'P2'}}, {\tt{'NED'}}, {\tt{'RWG'}} respectively, depending on the desired finite element in the preceding list.

Besides the definition of the finite element spaces, the toolbox {\sc{openFEM}} contains a few more functionalities, as the following.
\begin{itemize}
\item Operators. 
Operators can be specified on the finite element space itself. Available operators are:
\begin{itemize}
\item {\tt{grad, div, curl}}, which are differential operators that act on scalar or vectorial finite elements.
\item {\tt{curl, div, nxgrad, divnx, ntimes}}. Those operators are defined when solving problems on a bidimensional surface in $\mathbb{R}^3$. Here {\tt{n}} stands for the (outer) normal to the surface and all the differential operators are defined on surfaces. Such operators are commonly used when solving problems with the BEM (see below).  
\end{itemize}
\item Plots. Basic functions to plot a finite element or a solution are available. Namely, we have introduced
\begin{lstlisting}[style=Matlab-editor]
plot(Vh);
\end{lstlisting}
where {\tt{Vh}} is a finite element space. This produces a plot of the geometric location of the degrees of freedom that define the functions in {\tt{Vh}}.
\begin{lstlisting}[style=Matlab-editor]
surf(Vh,uh);
\end{lstlisting}
in order to plot a solution. In that case the figure produced consists in the geometry on which the finite element is defined colored by the magnitude of {\tt{uh}}. Eventually, as we have seen in the first examples presented in this paper,
the command {\tt{graph}} plots the graph of a finite element computed on a 2D flat surface. 
\begin{lstlisting}[style=Matlab-editor]
graph(Vh,uh);
\end{lstlisting}
\end{itemize}

\subsection{The {\tt{integral}} keyword ({\sc{openDOM}})}
Every integration done on a domain is evaluated through the keyword {\tt{integral}}. Depending on the context explained below, the returned value can be either a number, a vector or a matrix. More precisely, among the possibilities are
\begin{itemize}
\item {\tt{I = integral(dom, f);}}\\
where {\tt{dom}} is an integration domain on which the integral of the function {\tt{f}} needs to be computed. In that case the
function {\tt{f}} should be defined in a vectorial way, depending on a variable {\tt{X}} which can be a $N\times 3$ matrix.
For instance the definitions
\begin{lstlisting}[style=Matlab-editor]
f = @(X) X(:,1).*X(:,2);
g = @(X) X(:,1).^2 + X(:,2).^2 + X(:,3).^2;
h = @(X) X(:,1) + X(:,2).*X(:,3);
\end{lstlisting}
respectively stand for the (3 dimensional) functions
$$
f(x,y,z) = xy\,,\,\,\,g(x,y,z) = x^2+y^2+z^2\,,\,\,\,h(x,y,z) = x+yz\,.
$$
Since domains are all 3 dimensional (or more precisely embedded in the 3 dimensional space), only functions of 3 variables are allowed.
\item {\tt{I = integral(dom, f, Vh);}}\\
In that case, {\tt{f}} is still a 3 dimensional function as before while {\tt{Vh}} stands for a finite element space. The returned value {\tt{I}} is a column vector whose entries are given by 
$$
I_i =\int_{dom} f(X) \phi_i(X)\,dX
$$
for all basis function $\phi_i$ of the finite element space.
\item {\tt{I = integral(dom, Vh, f);}}\\
This case is identical to the previous one but now, the returned vector is a row vector.
\item {\tt{A = integral(dom, Vh, Wh);}}\\
where both {\tt{Vh}} and {\tt{Wh}} are finite element spaces. This returns the matrix $A$ defined by
$$
A_{ij} =\int_{dom} \phi_i(X) \psi_j(X)\,dX
$$
where $\phi_i$ (resp. $\psi_j$) stands for the current basis function of {\tt{Vh}} (resp. {\tt{Wh}}).
\item {\tt{A = integral(dom, Vh, f, Wh);}}\\
This is a simple variant where the entries of the matrix $A$ are now given by 
$$
A_{ij} =\int_{dom} f(X) \phi_i(X) \psi_j(X)\,dX\,.
$$
\end{itemize}

As a matter of fact, the leftmost occurring finite element space is assumed to correspond to test functions while the rightmost one corresponds to the unknowns.

\section{Generalization of the approach to the BEM}
It turns out that the preceding approach, described in section \ref{sec3},  can be generalized to the Boundary Element Method (BEM).
In such a context, after discretization, one has to solve a linear system where the underlying matrix is fully populated.
Typical examples are given by the acoustic or electromagnetic scattering. Indeed, let us consider a kernel\footnote{For instance, in the case of acoustic scattering in 3D, $G$ is the Helmholtz kernel defined by $G(x,y) = \frac{\exp(ik|x-y|)}{4\pi|x-y|}$.} $G(x,y)$
for which one has to compute the matrix $A$ defined by the entries
\begin{equation}
A_{ij} = \int_{x\in \Sigma_1}\int_{y\in \Sigma_2}\phi_i(x)G(x,y)\psi_j(y)\,dx\,dy\,
\end{equation}
the functions $(\phi_i)_{1\leq i\leq N_1}$ and  $(\psi_j)_{1\leq j\leq N_2}$, being basis functions of possibly different finite element spaces. Taking discrete integration formulas
on $\Sigma_1$ and $\Sigma_2$ respectively defined by the points and weights $(x_k,\omega_k)_{1\leq k\leq N_{int1}}$ and $(y_l,\eta_l)_{1\leq l\leq N_{int2}}$, leads to the approximation
\begin{equation}
A_{ij} \sim \sum_k\sum_l\phi_i(x_k)\omega_k G(x_k,y_l)\eta_l\psi_j(y_l)\,,
\end{equation}
which enables us to write in a matrix form
\begin{equation}
A\sim \Phi W_x G W_y \Psi\,.
\label{BEM}
\end{equation}
In this formula, the matrices $W_x$ and $W_y$ are the diagonal sparse matrices defined as before in (\ref{WB}) which depend on the quadrature weights $\omega_k$ and $\eta_l$ respectively. The matrices $\Phi$ and $\Psi$ are the (usually sparse) matrices defined as in (\ref{A=BWB}) for the basis functions $\phi_i$ and $\psi_j$ respectively, and $G$ is the dense matrix of size $N_{int1}\times N_{int2}$ given by $G_{kl} = G(x_k,y_l)$.

Again, building the sparse matrices as before, one only needs to further compute the dense matrix $G$ and assemble the preceding matrix $A$ with only matrix products.

\subsection{Generalization of the {\tt{integral}} keyword ({\sc{openDOM}})}
In terms of the syntax, we have extended the range of the {\tt{integral}} keyword in order to handle such integrals.
Indeed, the preceding formulas show that there are very little differences with respect to the preceding FEM
formulations. Namely, we now need to handle
\begin{itemize}
\item Integrations over 2 possibly different domains $\Sigma_1$ and $\Sigma_2$;
\item Any kernel depending on two variables provided by the user;
\item As before, two finite element spaces that are evaluated respectively on $x$ and $y$.
\end{itemize}
Furthermore, other formulations exist for the BEM, such as the so-called {\em{collocation method}}, in which one of the two integrals
is replaced by an evaluation at a given set of points. This case is also very much used when one computes (see the section below) the radiation of a computed solution on a given set of points.  

To handle all these situations three possible cases are provided to the user:
\begin{itemize}
\item The case with two integrations and two finite element spaces. This corresponds to computing the matrix
$$
A_{ij} = \int_{x\in \Sigma_1}\int_{y\in \Sigma_2} \phi_i(x)G(x,y)\psi_j(y)\,dx\,dy\,,
$$
and is simply performed in {\sc{Gypsilab}} by the following command.
\begin{lstlisting}[style=Matlab-editor]
A = integral(Sigma1, Sigma2, Phi, G, Psi);
\end{lstlisting}
As before, the first finite element space {\tt{Phi}} is considered as the test-function space while the second one, {\tt{Psi}}, stands for the unknown. The two domains on which the integrations are performed are given in the same order as the finite element spaces (i.e. {\tt{Phi}} and {\tt{Psi}} are respectively defined on $\Sigma_1$ and $\Sigma_2$).
\item The cases with only one integration and one finite element space. Two possibilities fall into this category. Namely, the computation of the matrix
$$
B_{ij} = \int_{x\in \Sigma_1} \phi_i(x)G(x,y_j)\,dx\,,
$$
for a collection of points $y_j$ and 
the computation of the matrix
$$
C_{ij} = \int_{y\in \Sigma_2} G(x_i,y)\psi_j(y)\,dy\,,
$$
for a collection of points $x_i$.
Both cases are respectively (and similarly) handled by the two following commands.\\
\begin{lstlisting}[style=Matlab-editor]
B = integral(Sigma1,y_j,Phi,G);
C = integral(x_i,Sigma2,G,Psi);
\end{lstlisting}
\end{itemize}
In all the preceding commands, $G$ should be a {\sc{Matlab}} function that takes as input a couple of 3 dimensional variables {\tt{X}} and {\tt{Y}} of respective sizes $N_X\times 3$ and $N_Y\times 3$. It should also be defined
in order to possibly handle sets of such points in a vectorized way. As a sake of example, $G(x,y) = \exp(i x\cdot y)$
can simply be declared as
\begin{lstlisting}[style=Matlab-editor]
G = @(X,Y) exp(1i*X*Y');
\end{lstlisting}
where it is expected that both $X$ and $Y$ are matrices that contain 3 columns (and a number of lines equal to the number of points $x$ and $y$ respectively).

\subsection{Regularization of the kernels}
It is commonly known that usual kernels that are used in classical BEM formulations (e.g. Helmholtz kernel in acoustics) are singular near 0. This creates a difficulty when one uses the BEM which may give very inaccurate results since the quadrature rules used for the $x$ and $y$ integration respectively may possess points that are very close one to another. However, the kernels often have a singularity which is asymptotically known.

In {\sc{Gypsilab}}, we provide the user with a way to regularize the considered kernel by computing a correction depending on its asymptotic behavior. As a sake of example, we consider the Helmholtz kernel which is used to solve the equations for acoustics (see section \ref{secacoustic} for much more details)
\begin{equation}
G(x,y)=\frac{e^{ik|x-y|}}{4\pi |x-y|}\,.
\label{Helmholtz}
\end{equation}
This kernel possesses a singularity when $x\sim y$ which has the asymptotics
$$
G(x,y)\sim_{x\sim y} \frac{1}{4\pi |x-y|} +O(1)\,.
$$
The idea is that the remainder is probably well approximated using Gauss quadrature rule, and we only need to correct the singular part coming from the integration of $\frac{1}{|x-y|}$. In {\sc{Gypsilab}}, this reads as
\begin{lstlisting}[style=Matlab-editor,basicstyle=\small]
A = 1/(4*pi)*(integral(Sigma,Sigma,Vh,Gxy,Vh));
A = A+1/(4*pi)*regularize(Sigma,Sigma,Vh,'[1/r]',Vh);
\end{lstlisting}
The first line, as we have already seen assembles the full matrix defined by the integral
$$
A_{ij}=\int_\Sigma \int_\Sigma G(x,y) \phi_i(x)\phi_j(y)\,dx\,dy\,, 
$$
where $(\phi_i)_i$ stands for the basis functions of the finite element space {\tt{Vh}}. 

The second line, however, regularizes the preceding integral by considering only the asymptotic behavior of $G$. 
This latter term computes and returns the {\bf{difference}} between an accurate computation and the Gauss quadrature evaluation of 
$$
\int_\Sigma \int_\Sigma \frac{\phi(x)\phi(y)}{4\pi|x-y|}\,dx\,dy\,.
$$
The Gauss quadrature is evaluated as before while the more accurate integration is computed using a semi-analytical method in which the integral in $y$ is computed analytically while the one in $x$ is done using a Gauss quadrature rule. The correction terms are only computed for pairs of integration points that are close enough. Therefore, the corresponding correction matrix is sparse.

\subsection{Coupling with {\sc{openHMX}}}
As it is well-known, and easily seen from the formula (\ref{BEM}), the matrices computed for the BEM are fully populated. Indeed, usual kernels $G(x,y)$ (e.g. the Helmholtz kernel) never vanish for any couple $(x,y)$. This
therefore leads to a matrix $G$ in (\ref{BEM}) for which no entry vanishes. Furthermore, the number of integration points $N_{int1}$ and
$N_{int2}$ is very often much larger than the corresponding numbers of degree of freedom. This means that the matrix $G$
usually has a size much larger than the final size of the (still fully populated) matrix $A$\footnote{As a sake of example, when one uses $P^1$ finite elements but an integration on triangles with 3 Gauss points per triangle, there are 6 times more Gauss points than unknowns (in a triangular mesh, the number of elements scales like twice the number of vertices). Calling $N$ the number of unknowns, the final matrix $A$ has a size $N^2$ while the matrix corresponding to the interaction of Gauss points is of size $(6N)^2=36 N^2$ which is much bigger.}. Both these facts limit very much the applicability of the preceding approach on classical computers to a number of degrees of freedom of a few thousand, which
is often not sufficient in practice. For this reason we also provide a coupling with the {\sc{Gypsilab}} toolbox {\sc{openHMX}} \cite{Gypsilab} in order to assemble directly a hierarchical $\mathcal{H}-$matrix compressed
version of the preceding matrices. Namely, for a given tolerance {\tt{tol}}, the commands

\begin{lstlisting}[style=Matlab-editor]
A = integral(Sigma1,Sigma2,Phi,G,Psi,tol);
B = integral(Sigma1,y_j,Phi,G,tol);
C = integral(x_i,Sigma2,G,Psi,tol);
\end{lstlisting}
return the same matrices as before, but now stored in a hierarchical $\mathcal{H}-$matrix format, and approximated to the desired tolerance. In particular, this enables the user to use the $+$, $-$, $*$, \textbackslash, {\tt{lu}} or {\tt{spy}} commands as if they were classical {\sc{Matlab}} matrix objects.  

These generalizations, together with the possibility of directly assemble $\mathcal{H}$-matrices using the same kind
of syntax, seem to us one of the cornerstones of the {\sc{Gypsilab}} package. To our knowledge, there is, at the moment, no comparable software package which handles BEM or compressed BEM matrices defined in a way as general and simple as here.

\subsection{Acoustic scattering}
\label{secacoustic}
As a matter of example, we provide hereafter the resolution of the acoustic scattering of a sound soft sphere 
and the corresponding program in {\sc{Gypsilab}}. For this test case, one considers a sphere of unit radius $\mathbb{S}^2$, and an incident acoustic wave given by
\begin{equation}
p_{inc}(x) = \exp(i k x\cdot d)
\end{equation}
where $k$ is the current wave number and $d$ is the direction of propagation of the wave. It is well-known that
the total pressure $p_{tot}$ outside the sphere is given by $p_{tot}=p_{inc}+p_{sca}$ where the scattered pressure wave obeys the formula
\begin{equation}
p_{sca}(x) = \int_{\mathbb{S}^2} G(x,y) \lambda(y)\,d\sigma(y)\,.
\end{equation}
In the preceding formula, the Green kernel of Helmholtz equation is given by
(\ref{Helmholtz}) and the density $\lambda$ is computed using the so-called single layer formula
\begin{equation}
- p_{inc}(x) = \int_{\mathbb{S}^2} G(x,y) \lambda(y)\,d\sigma(y)\,,
\label{Singlelayer}
\end{equation}
for $x\in \mathbb{S}^2$. This ensures that $p_{tot}=0$ on the sphere. Solving the equation (\ref{Singlelayer}) with the 
Galerkin finite element method amounts to solve the weak form
\begin{equation}
\int_{\mathbb{S}^2} \int_{\mathbb{S}^2} \mu(x) G(x,y) \lambda(y)\,d\sigma(x)\,d\sigma(y) = - \int_{\mathbb{S}^2} \mu(x) p_{inc}(x)\,d\sigma(x)\,,
\label{Singlelayerwf}
\end{equation}
where the test function $\mu$ and the unknown $\lambda$ belong to a discrete finite element space. We take the space $P^1$ defined on a triangulation $\mathcal{T}_h$ of $\mathbb{S}^2$. 
\vspace*{0.3cm}

\hrule
\begin{lstlisting}[style=Matlab-editor,basicstyle=\small]
% Parameters
N   = 1e3; 
tol = 1e-3;
X0  = [0 0 -1];
% Spherical mesh
sphere = mshSphere(N,1);
S2  = dom(sphere,3);    
% Radiative mesh - Vertical square
square = mshSquare(5*N,[5 5]);
square = swap(square);
% Frequency adjusted to maximum edge size
stp = sphere.stp;
k  = 1/stp(2)
f   = (k*340)/(2*pi)
% Incident wave
PW = @(X) exp(1i*k*X*X0');
% Green kernel: G(x,y) = exp(ik|x-y|)/|x-y| 
Gxy = @(X,Y) femGreenKernel(X,Y,'[exp(ikr)/r]',k);
% Finite element space
Vh = fem(sphere,'P1');
% Operator \int_Sx \int_Sy psi(x)' G(x,y) psi(y) dx dy 
LHS = 1/(4*pi)*(integral(S2,S2,Vh,Gxy,Vh,tol));
LHS = LHS+1/(4*pi)*regularize(S2,S2,Vh,'[1/r]',Vh);
% Wave trace --> \int_Sx psi(x)' pw(x) dx
RHS = integral(S2,Vh,PW);
% Solve linear system [-S] * lambda = - P0
lambda  = LHS \ RHS; 
% Radiative operator \int_Sy G(x,y) psi(y) dy 
Sdom = 1/(4*pi)*integral(square.vtx,S2,Gxy,Vh,tol);
Sdom = Sdom+1/(4*pi)*regularize(square.vtx,S2,'[1/r]',Vh);
% Domain solution : Pdom = Pinc + Psca
Pdom = PW(square.vtx) - Sdom * lambda;
% Graphical representation
figure
plot(square,abs(Pdom))
title('Total field solution')
colorbar
view(0,0);
hold off
\end{lstlisting}
\hrule
\vspace*{0.3cm}

\begin{figure}[h]
\centerline{\includegraphics[width = 8cm]{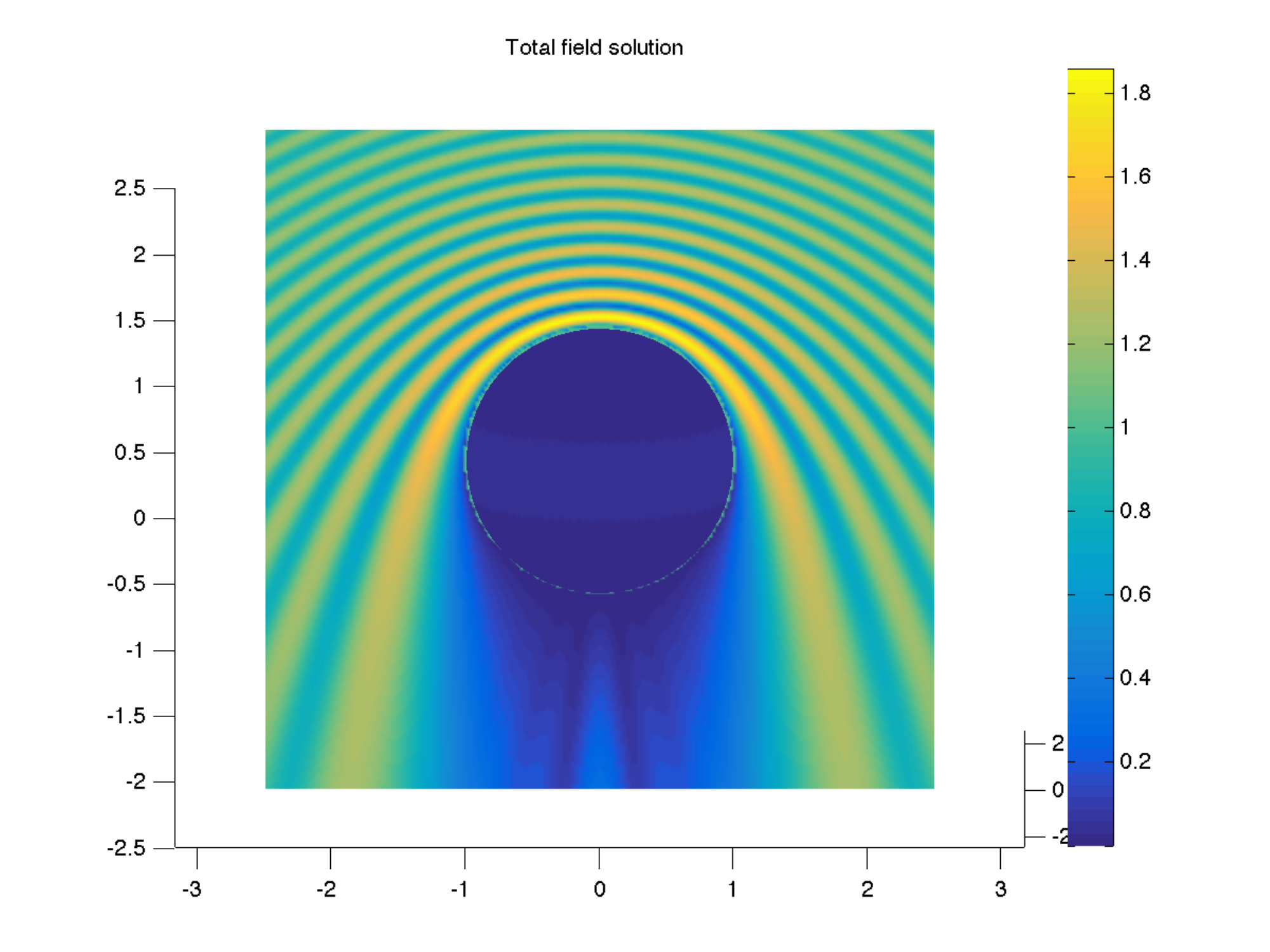}\includegraphics[width = 8cm]{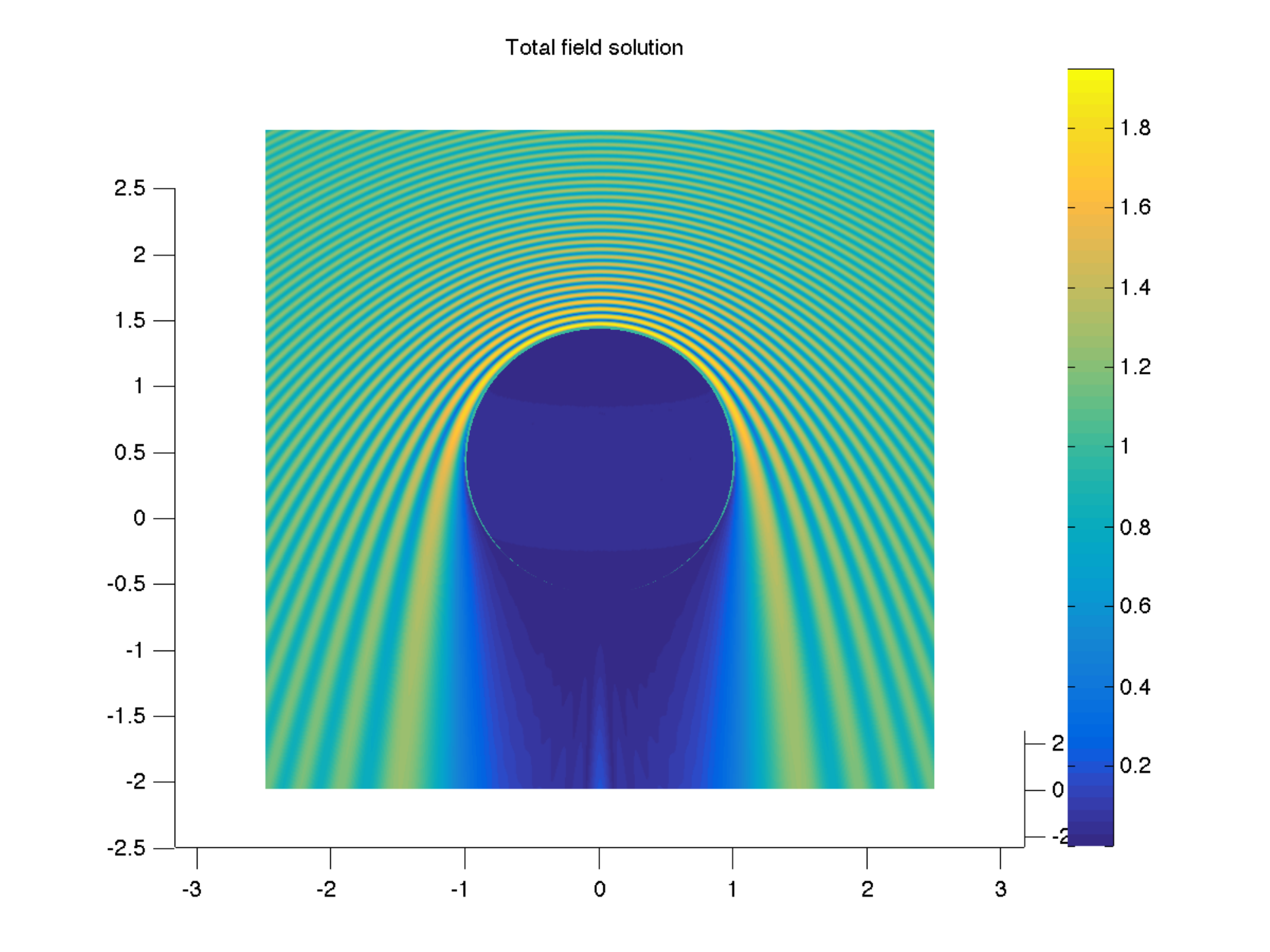}}
\caption{Magnitude of the pressure produced in the acoustic scattering by a unit sphere of a plane wave coming from above, using {\sc{Gypsilab}}. On the left, the sphere is discretized with $10^4$ vertices and the frequency is $10^3$ Hz. On the right the sphere is discretized with $9\cdot 10^4$ vertices and the frequency used is $3\cdot 10^3$ Hz.}
\label{fig3}
\end{figure}

\begin{figure}[h]
\centerline{\includegraphics[width = 9cm]{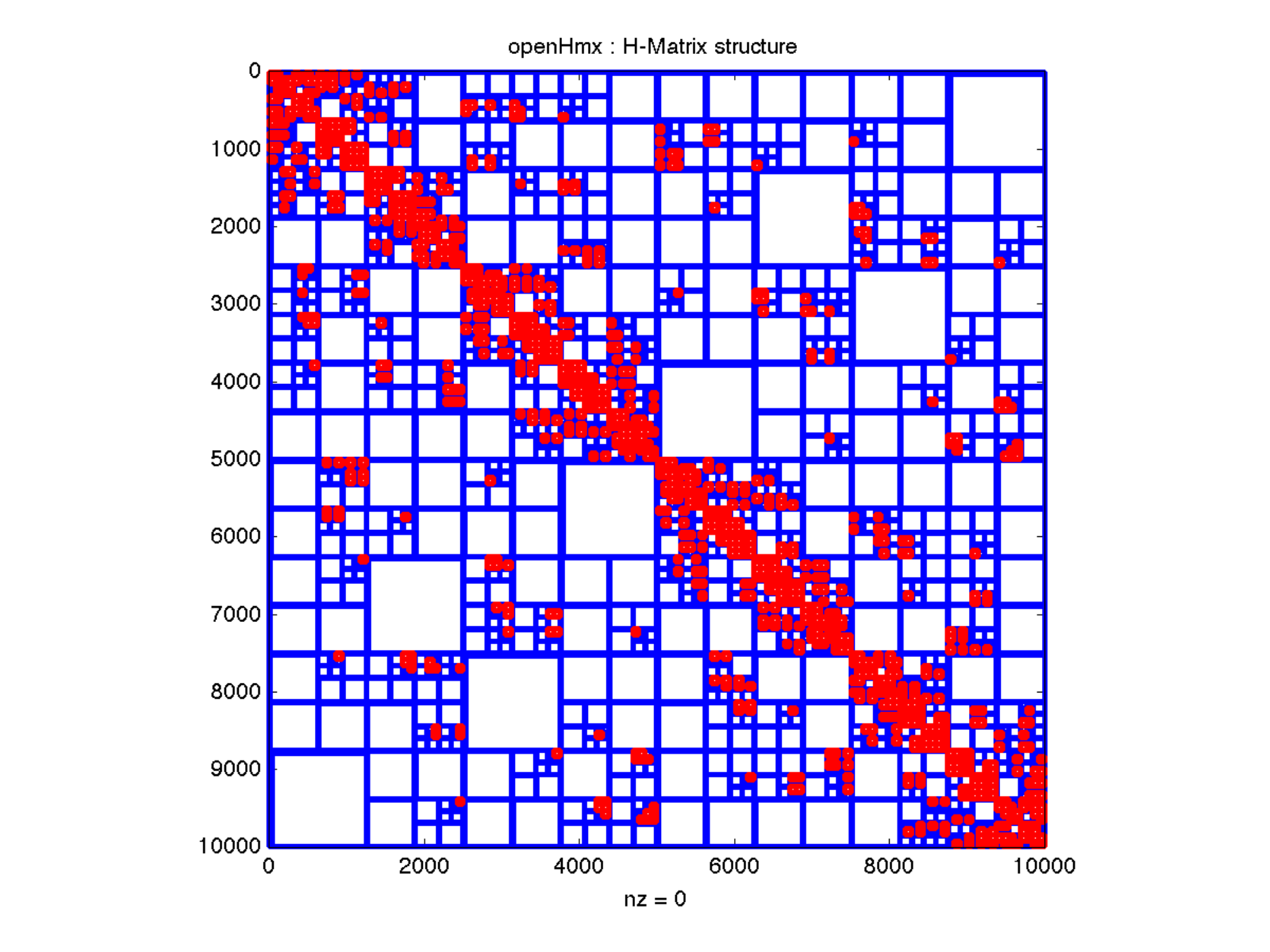}\hspace*{-1cm}\includegraphics[width = 9cm]{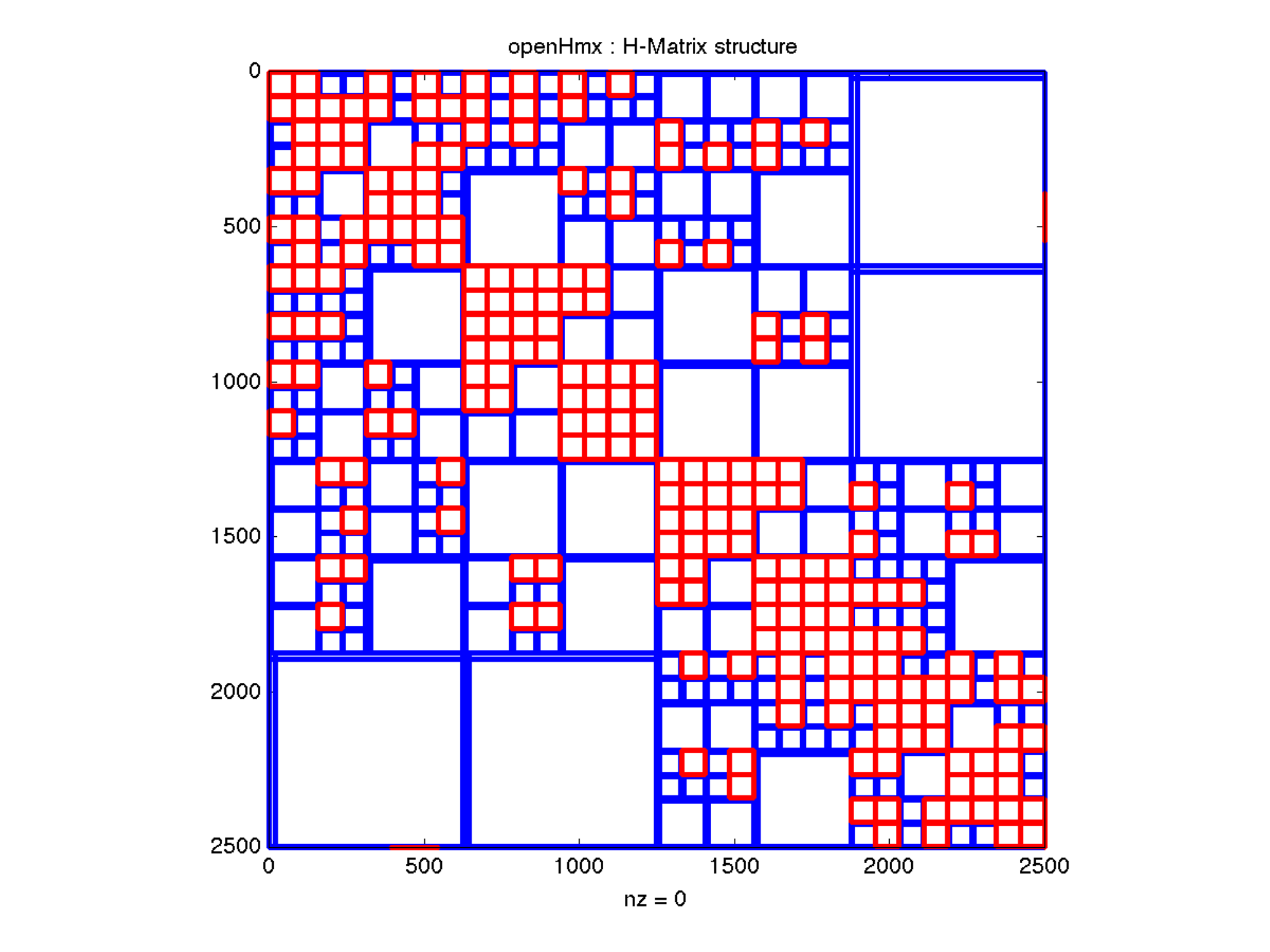}}
\caption{The $\mathcal{H}$-matrix produced in the case of the acoustic scattering with $10^4$ unknowns. The left-hand side picture is obtained by using the {\tt{spy}} command on the matrix itself. A zoom on the upper left part of the matrix (right) shows that each block contains an information about the local rank.}
\label{Hmat}
\end{figure}

The preceding program follows the traditional steps for solving the problem.
Namely, one recognizes the spherical mesh and domain (lines 11-12), the radiative mesh on which we want to compute and plot the solution, here a square (lines 14-15),
the incident plane wave (line 21), the Green kernel definition (line 23), the finite element space (line 25), the assembling 
of the operator (line 27-28), the construction of the right-hand side (line 30), and the resolution
of the problem (line 32). The rest of the program consists in computing from the solution $\lambda$ of (\ref{Singlelayer}), the
total pressure on the radiative mesh, and plot it on the square mesh. 
Notice that due to the presence of the {\tt{tol}} parameter in the assembling of the operator
(and also of the radiative operator), the corresponding matrices are stored as $\mathcal{H}-$matrices. Notice also that  the key part of the method (assembling and resolution) are completely contained between lines 21-32. The figures of the total pressure are given in Fig. \ref{fig3} for the two cases $N=10^4$ and $N=9\cdot 10^4$. The $\mathcal{H}$-matrix produced in the former case is shown in Fig. \ref{Hmat}.

\subsection{Electromagnetic scattering}
In electromagnetic scattering, formulations involving integral equations discretized with the BEM are also commonly used. We refer the reader to \cite{Colton,Nedelec} for an overview of classical properties of integral operators and discretizations. Three formulations are used to compute the magnetic current $J=n\times H$ on the surface of the scatterer. Namely, we distinguish
\begin{itemize}
\item The Electric Field Integral Equation (EFIE)
$$
TJ = -E_{inc,t}
$$
where the single layer operator $T$ is defined by
$$
TJ = ik \int_{\Sigma} G(x,y) J(y) \,dy + \frac{i}{k}\nabla_x \int_{\Sigma} G(x,y) \mbox{div} J(y) \,dy
$$
and $E_{inc,t}$ is the tangential component of the incident electric field.
\item The Magnetic Field Integral Equation (MFIE)
$$
\left(\frac12 - n\times K\right) J = -n\times H_{inc,t}
$$
where the double layer operator $K$ is defined by
$$
KJ = \int_{\Sigma} \nabla_y G(x,y) J(y) \,dy
$$
and $H_{inc,t}$ is the tangential component of the incident magnetic field.
\item The Combined Field Integral Equation (CFIE), used to prevent the ill-posedness of the preceding formulations
at some frequencies. It is a linear combination of the Electric and Magnetic Field Integral Equations
$$
\left(-\beta T +(1-\beta)\left(\frac12 - n\times K\right) \right)J = \beta E_{inc,t} -(1-\beta)n \times H_{inc,t}\,.
$$
\end{itemize}
As before, the kernel is the Helmholtz Green kernel defined by (\ref{Helmholtz}).

The classical finite element formulation for this problem uses the Raviart-Thomas elements for $J$ which are available in {\sc{openFEM}}.
The key part of the program assembling the CFIE operator and solving the scattering problem is given hereafter. For simplicity, we only focus on the assembling and solving parts and do not provide the initialization part or the radiation and plotting parts.

\vspace*{0.3cm}

\hrule
\begin{lstlisting}[style=Matlab-editor,basicstyle=\small]
% Incident direction and field
X0 = [0 0 -1]; 
E  = [0 1  0];    % Polarization of the electric field
H  = cross(X0,E); % Polarization of the magnetic field
% Incident Plane wave (electromagnetic field)
PWE{1} = @(X) exp(1i*k*X*X0') * E(1);
PWE{2} = @(X) exp(1i*k*X*X0') * E(2);
PWE{3} = @(X) exp(1i*k*X*X0') * E(3);
%
PWH{1} = @(X) exp(1i*k*X*X0') * H(1);
PWH{2} = @(X) exp(1i*k*X*X0') * H(2);
PWH{3} = @(X) exp(1i*k*X*X0') * H(3);
% Green kernel function G(x,y) = exp(ik|x-y|)/|x-y| 
Gxy    = @(X,Y) femGreenKernel(X,Y,'[exp(ikr)/r]',k);
Hxy{1} = @(X,Y) femGreenKernel(X,Y,'grady[exp(ikr)/r]1',k) ;
Hxy{2} = @(X,Y) femGreenKernel(X,Y,'grady[exp(ikr)/r]2',k) ;
Hxy{3} = @(X,Y) femGreenKernel(X,Y,'grady[exp(ikr)/r]3',k) ;
% Finite elements
Vh = fem(sphere,'RWG');
% Finite element mass matrix
Id = integral(sigma,Vh,Vh);
% Finite element boundary operator
T = 1i*k/(4*pi)*integral(sigma,sigma,Vh,Gxy,Vh,tol) ...
 -1i/(4*pi*k)*integral(sigma,sigma,div(Vh),Gxy,div(Vh),tol);
T = T + 1i*k/(4*pi)*regularize(sigma,sigma,Vh,'[1/r]',Vh) ...
 -1i/(4*pi*k)*regularize(sigma,sigma,div(Vh),'[1/r]',div(Vh));
% Finite element boundary operator
nxK = 1/(4*pi)*integral(sigma,sigma,nx(Vh),Hxy,Vh,tol); 
nxK = nxK+1/(4*pi)*regularize(sigma,sigma,nx(Vh),'grady[1/r]',Vh);
% Left hand side
LHS = -beta * T + (1-beta) * (0.5*Id - nxK);
% Right hand side
RHS = beta*integral(sigma,Vh,PWE) ...
        - (1-beta)*integral(sigma,nx(Vh),PWH);
% Solve linear system 
J = LHS \ RHS;
\end{lstlisting}
\hrule
\vspace*{0.3cm}

As one can see, the program is a direct transcription of the mathematical weak formulation of the problem. This follows the same lines as in the acoustic scattering problem except for the operators that are different and the finite element used. Notice also the regularization of the double layer kernel in line 29.

\section{Performances}
\subsection{Performances in FEM}
In this section we compare {\sc{Gypsilab}} with {\sc{FreeFem++}}. The machine that we have used for this comparison is equipped with Xeon E5-2643-V1 processors with a frequency of 3.3 GHz and 128 GB of memory. Although the machine possesses two such processors, meaning that up to 8 cores could be used for the computation, we only chose a single core to run the test, both for {\sc{FreeFem++}} and {\sc{Matlab}} which is therefore launched using the {\tt{-singleCompThread}} option. For the test, we have used {\sc{FreeFem++}} version 3.61 and {\sc{Matlab}} R2018a. 

We have chosen to solve the Dirichlet problem
\begin{equation}
\left\{
\begin{array}{l}
-\Delta u = 1\,\mbox{ in } \Omega \,,\\
u = 0 \,\mbox{ on } \partial \Omega\,,
\end{array}
\right.
\label{laplacetestp12}
\end{equation}
where $\Omega=(0,1)^3$ is the unit cube of $\mathbb{R}^3$ meshed regularly, with $(N+1)^3$ points and $N$ ranges from 20 up to 100 depending on the case. Both linear $P^1$ and quadratic $P^2$ elements are considered.
In all cases, the same quadrature formula was used to approximate the underlying integrals, namely with one (resp. 4) quadrature point(s) for $P^1$ (resp. $P^2$) elements. Tables \ref{table1} and \ref{table2} give the respective timings for 
assembling the problem and solving it using in both cases the GMRES solver and a $10^{-6}$ accuracy for the convergence.
Notice also that the number of degrees of freedom $N_{dof}$ is different in {\sc{FreeFem++}} and {\sc{Gypsilab}}. This is due to the fact that
Dirichlet boundary conditions are enforced by penalization in {\sc{FreeFem++}}, while we eliminate the corresponding unknowns in 
{\sc{Gypsilab}}.

We notice that {\sc{Gypsilab}} appears slower than {\sc{FreeFem+}} by a factor which is less than 4 in all the configurations
that were tested. We therefore believe that {\sc{Gypsilab}} is a suitable tool for prototyping.
 
\begin{table}[h]
\begin{center}
\begin{tabular}{c||c|c|c||c|c|c}
$P^1$ & \multicolumn{3}{c||}{{\sc{FreeFem++}}} & \multicolumn{3}{c}{{\sc{Gypsilab}}}\\ 
\hline
$N$ & $N_{dof}$ & $T_{ass}$ & $T_{sol}$ & $N_{dof}$ & $T_{ass}$ & $T_{sol}$  \\
\hline
\hline
20 &  9 261 & 0.28 & 0.10 & 6859 & 0.36 & 0.15 \\
\hline
40 & 68 921 & 0.34 & 0.86 & 59319 & 2.3 & 1.7 \\
\hline
60 & 226 981 & 1.1 & 4.3 & 205379 & 8.0 & 10.8\\
\hline
80 & 531 441 & 2.5 & 20.3 & 493039 & 20 & 49\\
\hline
100 & 1 030 301 & 5.3 & 64.8 & 970299 & 43 & 230\\
\hline
\end{tabular}
\end{center}
\caption{Timings for assembling the matrix and solving the linear system coming from the discretization of the Laplace problem (\ref{laplacetestp12}) with $P^1$ finite elements in {\sc{FreeFem++}} and {\sc{Gypsilab}} respectivley. We notice a much bigger time for assembling the problem which is compensated by the resolution. {\sc{Gypsilab}} appears slower than {\sc{FreeFem++}} for the total resolution by a factor less than 4 in all cases.}
\label{table1}
\end{table}

\begin{table}[h]
\begin{center}
\begin{tabular}{c||c|c|c||c|c|c}
$P^2$ & \multicolumn{3}{c||}{{\sc{FreeFem++}}} & \multicolumn{3}{c}{{\sc{Gypsilab}}}\\ 
\hline
$N$ & $N_{dof}$ & $T_{ass}$ & $T_{sol}$ & $N_{dof}$ & $T_{ass}$ & $T_{sol}$  \\
\hline
\hline
20 &  68921 & 0.5 & 1.0 & 59319 & 3.0 & 1.8 \\
\hline
40 & 531441 & 3.8 & 28.36 & 493039 & 25.7 & 50.8 \\
\hline
60 & 1771561 & 13.0 & 202.0 & 1685159 & 91.7 & 507\\
\hline
\end{tabular}
\end{center}
\caption{Timings for assembling the matrix and solving the linear system coming from the discretization of the Laplace problem (\ref{laplacetestp12}) with $P^2$ finite elements in {\sc{FreeFem++}} and {\sc{Gypsilab}} respectivley. As before assembling the problem is much slower in {\sc{Gypsilab}}. However the total resolution times are comparable.}
\label{table2}
\end{table}

\subsection{Performances in BEM}
We report in this section the performances attained by the acoustic scattering of the sphere previously described. Here, the goal is not to compare with another package, but we have still used a single core iof the machine to run the test.
We give hereafter the timings for different meshes of the sphere corresponding to an increasing number $N$ 
of degrees of freedom in the underlying system. The first part of Table~ \ref{table4} gives the timings to assemble the full BEM matrix for sizes ranging from 1000 to 150000 degrees of freedom. Above 10000 unknowns, the matrix of the kernel computed
at the integration points no longer fits into the memory of the available machine. Therefore, we turn to use the hierarchical compression for the matrix, i.e. the $\mathcal{H}-$matrix paradigm available through the use of {\sc{openHMX}}. This enables us to increase the size of reachable problems by an order of magnitude and more. This is reported in the bottom part of Table \ref{table4}. Notice that the frequency of the problem is adapted to the precision of the mesh as shown in the last column of the table. 

In order to see the effect of the frequency on the construction of the $\mathcal{H}$-matrix 
and the resolution, we have also tried to fix the frequency to 316 Hz (the smallest value for the preceding case) and
check the influence on the assembling, regularization and solve timings in the problem. The data are given in Table \ref{table5}.
It can be seen that the underlying matrix is much easier compressed and quicker assembled. The resolution time is also significantly reduced. Indeed, the time to assemble the $\mathcal{H}$-matrix becomes proportional to the number of unknowns.

\begin{table}[ht]
\begin{center}
\begin{tabular}{c|c|c|c||c|c|c|c}
\hline
$N_{dof}$ & $T_{ass}$ & $T_{reg}$ & $T_{sol}$ & $T^H_{ass}$ & $T^H_{reg}$ & $T^H_{sol}$ & Freq. (Hz)\\
\hline
1000 & 3.43 & 1.46 & 0.19 & 5.58 & 1.49 & 0.69 & 316\\
3000 & 27.1 & 2.82 & 2.32 & 21.0 & 4.20 & 3.12 & 547\\
5000 & 74.0 & 4.72 & 8.20 & 31.6 & 5.13  & 6.54 & 707 \\
10000 & 318 & 9.80 & 64.0 & 72.9 & 11.3 & 20.1 & 1000 \\
20000 & -- & -- & -- & 163 & 22.8 & 70.2 & 1414\\
40000 & -- & -- & -- & 358 & 41.3 & 298 & 2000\\
80000 & -- & -- & -- & 1035 & 93.0 & 1230 & 2828 \\
150000 & -- & -- & -- &  3400 & 167 & 4728 & 3872\\
\hline
\end{tabular}
\end{center}
\caption{Timings in seconds for assembling, regularizing and solving the problem of acoustic scattering given in section \ref{secacoustic}. The second half of the table corresponds to the timings using the $\mathcal{H}-$matrix approach using in the {\sc{Gypsilab-openHMX}} package. For problems of moderate size it is slightly faster to use the classical BEM approach, while the sizes corresponding to the bottom lines are beyond reach for this method. Notice that when we solve this problem using the $\mathcal{H}-$matrices, a complete LU factorization is computed. This is not optimal in the present case, in particular when compared to other compression techniques such as the FMM, since the underlying linear system is solved only once (with only one right hand side). Non available data indicates that the problem cannot fit into memory.} 
\label{table4}
\end{table}

\begin{table}[ht]
\begin{center}
\begin{tabular}{c|c|c|c|c}
\hline
$N_{dof}$ & $T^H_{ass}$ & $T^H_{reg}$ & $T^H_{sol}$ & Freq. (Hz)\\
\hline
10000 & 60.1 & 9.32 & 10.0 & 316 \\
20000 & 118 & 21 & 22.3 & 316\\
40000 & 181 & 38.9 & 53.0 & 316\\
80000 & 389 & 79.3 & 133 & 316 \\
150000 & 672 & 141 & 279 & 316\\
\hline
\end{tabular}
\end{center}
\caption{Timings in seconds for assembling, regularize and solve the problem of acoustic scattering given in section \ref{secacoustic} at a fixed frequency $f=316 Hz$ on the unit sphere. The number of unknowns is given as $N_{dof}$ and the $\mathcal{H}$-matrix compression technique is used.} 
\label{table5}
\end{table}
\section{Conclusion}
The package {\sc{Gypsilab}} is a numerical library written in full {\sc{Matlab}} that allows the user to solve PDE using the finite element technique.
Very much inspired by the {\sc{FreeFem++}} formalism, the package contains classical FEM and BEM functionalities. In this latter case, the library allows the user to store the operators in a $\mathcal{H}$-matrix format that makes it possible to factorize the underlying matrix and solve using a direct method the linear system. We are not aware of any comparable software package that combines ease of use, generality and performances to the level reached by {\sc{Gypsilab}}. We have shown illustrative examples in several problems ranging from classical academic Laplace problems to the Combined Field Integral Equation in electromagnetism. Eventually, a short performance analysis shows that the library possesses enough performance to run problems with a number of unknowns of the order of a million in reasonable times. A lot remains to be done, as extending the available finite elements, proposing different compression strategies or coupling FEM and BEM problems, that we wish to study now. In particular solving coupled FEM-BEM problems in {\sc{Gypsilab}} will be the subject of a forthcoming paper.

Finally, {\sc{Gypsilab}} is available under GPL license \cite{Gypsilab} and which makes it a desirable tool for prototyping. 

\section*{Acknowledgments} 
Both authors would like to thank Pierre Jolivet for valuable comments on the present work, especially concerning the
way {\sc{FreeFem++}} handles different quadrature formulas. The financial support of the french Direction G\'en\'erale de l'Armement is also gratefully acknowledged.

\end{document}